\def\xi{{\xi}}
\def\p{{\pi}}

\def\Ga{{\Gamma}}

\def\Xi{{\Xi}}
\def\Pi{{\Pi}}

\def\B{{\cal B}}

\def\aF{{\mathbb F}}
\def\aZ{{\mathbb Z}}
\def\bB{{\bm B}}
\def\bd{{\bm d}}
\def\be{{\bm e}}
\def\bF{{\bm F}}

\def\bmf{{\bm f}}
\def\bg{{\bm g}}
\def\bn{{\bm n}}
\def\bp{{\bm p}}
\def\bu{{\bm u}}
\def\bv{{\bm v}}
\def\bw{{\bm w}}
\def\bx{{\bm\chi}}

\def\fn{{\scriptstyle}}
\def\ube{{\fn\be}}
\def\ubmf{{\fn\bmf}}
\def\ubF{{\fn\bF}}
\def\ubn{{\fn\bn}}
\def\ubu{{\fn\bu}}
\def\ubv{{\fn\bv}}
\def\ux{{\fn\bx}}

\def\emt{{\emptyset}}
\def\op{{\oplus}}

\def\pf{{\hfill$\Box$}}

\def\la{{\langle}}
\def\ra{{\rangle}}
\def\c{{\cdot}}

\def\mod{{\ \rm mod\ }}

\def\ord{{\rm Ord}}
\def\set{{\rm Set}}
\def\val{{\rm\bf Val}}

\documentclass[12pt]{article}

\usepackage{psfig}
\usepackage{amssymb}

\newcommand{\bm}[1]{{\mbox{\boldmath $#1$}}}

\newtheorem{thm}{Theorem}

\newtheorem{lem}[thm]{Lemma}

\newtheorem{res}[thm]{Result}

\begin{document}

\baselineskip=24pt

%
%
\title{An Application of Graph Pebbling\\ to Zero-Sum Sequences in 
Abelian Groups}
\author{
 \\
Shawn Elledge\\
Department of Mathematics and Statistics\\
Arizona State University\\
Tempe, AZ 85287-1804\\
email: sme13@asu.edu\\
 \\
and \\
 \\
Glenn H. Hurlbert\thanks{Partially supported by National Security
    grant \#MDA9040210095.}\\
Department of Mathematics and Statistics\\
Arizona State University\\
Tempe, AZ 85287-1804\\
email: hurlbert@asu.edu\\
}
\maketitle
\newpage


%
%
\begin{abstract}
A sequence of elements of a finite group $G$ is called a zero-sum sequence 
if it sums to the identity of $G$.
The study of zero-sum sequences has a long history with many important
applications in number theory and group theory.
In 1989 Kleitman and Lemke, and independently Chung, proved a strengthening
of a number theoretic conjecture of Erd\H{o}s and Lemke.
Kleitman and Lemke then made more general conjectures for finite groups,
strengthening the requirements of zero-sum sequences.
In this paper we prove their conjecture in the case of abelian groups.
Namely, we use graph pebbling to prove that for every sequence 
$(g_k)_{k=1}^{|G|}$ of $|G|$ elements of a finite abelian group $G$ 
there is a nonempty subsequence $(g_k)_{k\in K}$ such that 
$\sum_{k\in K}g_k=0_G$ and $\sum_{k\in K}1/|g_k|\le 1$, where $|g|$ is 
the order of the element $g\in G$.
\vspace{0.2 in}

\noindent {\bf 2000 AMS Subject Classification:} 
11B75, 20K01, 05D05
\vspace{0.2 in}

\noindent {\bf Key words:} 
Graph pebbling, finite abelian group, zero-sum sequence
\end{abstract}
\newpage

%
%
\section{Introduction}\label{intro}

A sequence of elements of a finite group $G$ is called a zero-sum sequence 
if it sums to the identity of $G$.
A standard pigeonhole principle argument shows that any sequence of $|G|$
elements of $G$ contains a zero-sum subsequence; in fact having consecutive
terms (one can instead stipulate that the zero-sum subsequence has at most
$N$ terms --- where $N=N(G)$ is the exponent of $G$, i.e. the maximum order 
of an element of $G$ --- which is best possible).

First considered in 1956 by Erd\H{o}s \cite{Erd},
the study of zero-sum sequences has a long history with many important
applications in number theory and group theory.
In 1961 Erd\H{o}s et al. \cite{EGZ} proved that every sequence of $2|G|-1$
elements of a cyclic group $G$ contains a zero-sum subsequence of length
exactly $|G|$.
In 1969 van Emde Boas and Kruyswijk \cite{EBK} proved that any sequence of 
$N(1+\log(|G|/N))$ elements of a finite abelian group contains a zero-sum 
sequence.
In 1994 Alford et al. \cite{AGP} used this result and modified Erd\H{o}s's
arguments to prove that there are infinitely many Carmichael numbers.
Much of the recent study has involved finding Davenport's constant $D(G)$,
defined to be the smallest $D$ such that every sequence of $D$
elements contains a zero-sum subsequence \cite{O}.
Applications of the wealth of results on this problem 
\cite{Car,Gao,GG,GT,GS,Sun} and its variations \cite{GJ,Nat}
to factorization theory and to graph theory can be found in \cite{AFK,Cha}.

In 1989 Kleitman and Lemke \cite{KL}, and independently Chung \cite{Chu}, 
proved the following strengthening of a number theoretic conjecture of 
Erd\H{o}s and Lemke (see also \cite{CHH,D}).

\begin{res}\label{cyclic}
For any positive integer $n$, every sequence $(a_k)_{k=1}^n$ of $n$ 
integers contains a nonempty subsequence $(a_k)_{k\in K}$ such that 
$\sum_{k\in K}a_k\equiv 0\mod n$ and $\sum_{k\in K}\gcd(a_k,n)\le n$.
\end{res}

Kleitman and Lemke then made more general conjectures for finite groups,
strengthening the requirements of zero-sum sequences.
In this paper we prove their conjecture in the case of abelian groups.
Namely, we use graph pebbling (and Result \ref{cyclic}) to prove the
following theorem (we use $|g|$ to denote the order of the element
$g\in G$).

\begin{thm}\label{groups}
For every sequence $(g_k)_{k=1}^{|G|}$ of $|G|$ elements of a finite 
abelian group $G$ there is a nonempty subsequence $(g_k)_{k\in K}$ such 
that $\sum_{k\in K}g_k=0_G$ and $\sum_{k\in K}1/|g_k|\le 1$.
\end{thm}

Notice that Result \ref{cyclic} is the special case of Theorem \ref{groups}
in which $G$ is cyclic.
Also notice that the condition on the sum of the orders implies that
$|K|\le N(G)$, with equality if and only if $|g_k|=N$ for every $k\in K$.

%
%
\section{Preliminaries}\label{prelim}

%
%
\subsection{Graph Pebbling}\label{pebbling}

Let $\Ga=(V,E)$ be a graph with vertices $V$ and edges (unordered pairs
of edges) $E$.
Given a configuration of pebbles on $V$, a pebbling step consists of
removing two pebbles from a vertex $u$ and placing one pebble on an
adjacent vertex $v$ ($uv\in E$).
The pebbling number $\p=\p(\Ga)$ is the smallest number $\p$ such that, 
from every configuration of $\p$ pebbles on $V$ it is possible to place a
pebble on any specified target vertex after a sequence of pebbling moves.
There is a rapidly growing literature on graph pebbling 
\cite{CH1,CHKT,Her}, including variations such as optimal pebbling 
\cite{FW,Moe,PSV}, pebbling thresholds \cite{BBCH,BH,CH2} and cover pebbbling 
\cite{CCFHPST,HM,WY}.

One variation of graph pebbling involves labelling the edges $uv\in E$ 
by positive integer weights $w(uv)$, so that a pebbling step from $u$ to 
$v$ removes $w(uv)$ pebbles from $u$ before placing one pebble on $v$.
In this light, standard graph pebbling has weight 2 on every edge.
Let $\B^n$ be the graph of the $n$ dimensional boolean algebra ---
its vertices are all binary $n$-tuples; its edges are the pairs of
$n$-tuples that differ by a single digit.
For every edge between vertices that differ in the $i^{\rm th}$ digit,
let $w_i$ be its weight.
Finally, let $\bw=\la w_i\ra_{i=1}^n$ and denote the resulting weighted graph
by $\B^n(\bw)$.
Then Chung's theorem \cite{Chu} is as follows.

\begin{thm}
The generalized pebbling number of the weighted graph $\B^n(\bw)$ is
$\p(\B^n(\bw))=\prod_{i=1}^nw_i$.
\end{thm}

%
%
\subsection{Group Structure}\label{group}

Let $\aZ_n$ denote the finite cyclic group on $n$ elements.
The standard representation for an abelian group $G$ has the form
$\aZ_{N_1}\op\aZ_{N_2}\op\ldots\op\aZ_{N_r}$, where
$N_i|N_{i-1}$ for $1<i\le r$ (although, purposely, we've written
the order of the cycles in reverse to the standard).
Thus the exponent of $G$ is $N(G)=N_1$ and the rank of $G$ is $r(G)=r$.
One of the useful techniques in this paper is to break down each cycle
$\aZ_{N_i}$ into products of cycles of distinct prime powers.
We write $G=\op_{i=1}^t\op_{j=1}^{m_i}\aZ_{p_i^{e_{i,j}}}$ for
some primes $p_i$, multiplicities $m_i$, and exponents $e_{i,j}$.
Thus $G$ can be coordinatized so that elements $g$ have the form
$\bg=\la g_{i,j}\ra$, and addition is coordinatewise with the
$(i,j)^{\rm th}$ coordinate computed modulo $p_i^{e_{i,j}}$.
Further, instead of writing the primes $p_i$ in increasing order, we 
write them so that $e_{i,1}\ge\cdots\ge e_{i,m_i}$ for every $1\le i\le t$.
Hence the exponent of $G$ can be written
$N=N(G)=\prod_{i=1}^t p_i^{e_{i,1}}$.

%
%
\subsection{Notation}\label{notation}

\begin{figure}
\centerline{\psfig{figure=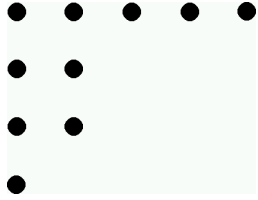,height=1.0in}}
\caption{Ferrer's diagram for $(5,2,2,1)$}\label{ferrer}
\end{figure}

As already witnessed, we will adopt the convention that bold fonts will 
denote vectors.
Let $\be_i=\la e_{i,j}\ra_{j=1}^{m_i}$, $\be=\la\be_i\ra_{i=1}^t$
and $m=\sum_{i=1}^t m_i$.
Then $\be_i$ can be thought of as a partition of the exponent of
$p_i$ in the prime factorization of $|G|$.
Define $\bd_i$ to be the dual partition that arises from the
Ferrer's diagram of $\be_i$.
For example, Figure \ref{ferrer} shows the Ferrer's diagram for $(5,2,2,1)$ 
(dots per row) and its dual $(4,3,1,1,1)$ (dots per column), both partitions 
of 10.

Next define $\bmf_{i,r}=\la {\bm 1}^r,{\bm 0}^{m_i-r}\ra$, and let
$$\bF_{i,r}
=\la\bmf_{1,0},\cdots,\bmf_{i-1,0},\bmf_{i,r},\bmf_{i+1,0},\cdots,\bmf_{m,0}\ra
=\la {\bm 0}^a,\bmf_{i,r},{\bm 0}^b\ra,\ $$ 
where $a=\sum_{i<r}m_i$ and $b=\sum_{i>r}m_i$.
For vectors $\bu=\la u_k\ra_{k=1}^s$, $\bv=\la v_k\ra_{k=1}^s$ and
$\bw=\la w_k\ra_{k=1}^s$ denote $\bu\ubv=\la u_k^{v_k}\ra_{k=1}^s$
and $\bu^{\c\ubv}=\prod_{k=1}^s u_k^{v_k}$.
Now let $\bp_i=\la p_i\ra_{j=1}^{m_i}$, $\bp=\la\bp_i\ra_{i=1}^t$ and
$\bp_0=\la p_i\ra_{i=1}^t$, and
define $\bn=\la n_i\ra_{i=1}^t=\la e_{i,1}\ra_{i=1}^t$ and
$n=\sum_{i=1}^t n_i$.
Note that $\bp_0^{\c\ubn}=N(G)$ and $\bp^{\c\ube}=|G|$.
We also write $\bu\le\bv$ when $u_k\le v_k$ for every $k$,
$\bu\equiv\bv\mod\bw$ when $u_k\equiv v_k\mod w_k$ for every $k$,
and
$\bu\bv=\bw$ (or $\bu=\bw/\bv$) when $u_kv_k=w_k$ for every $k$.

\begin{figure}
\begin{center}
$\be=\la 5,4,3,1;2,2;3;4,1,1\ra\ ,\quad
\be_1=\la 5,4,3,1\ra\ ,\quad
\bd_1=\la 4,3,3,2,1\ra$
\begin{tabular}{lclcl}
\hline
&&&&\\
$\be(0,0,0,0)$&$=$&&&$\la 5,4,3,1;2,2;3;4,1,1\ra$\\
&&&&\\
$\be(1,0,0,0)$&$=$&$\be(0,0,0,0)-\bF_{1,d_{1,u_1}}$&$=$&
	$\la 4,3,2,0;2,2;3;4,1,1\ra$\\
&&&&\\
$\be(1,1,0,0)$&$=$&$\be(1,0,0,0)-\bF_{2,d_{2,u_2}}$&$=$&
	$\la 4,3,2,0;1,1;3;4,1,1\ra$\\
&&&&\\
$\be(1,1,0,1)$&$=$&$\be(1,1,0,0)-\bF_{4,d_{4,u_4}}$&$=$&
	$\la 4,3,2,0;1,1;3;3,0,0\ra$\\
&&&&\\
$\be(2,1,0,1)$&$=$&$\be(1,1,0,1)-\bF_{1,d_{1,u_1}}$&$=$&
	$\la 3,2,1,0;1,1;3;3,0,0\ra$\\
&&&&\\
$\be(3,1,0,1)$&$=$&$\be(2,1,0,1)-\bF_{1,d_{1,u_1}}$&$=$&
	$\la 2,1,0,0;1,1;3;3,0,0\ra$\\
&$\cdot$&&$\cdot$&\\
&$\cdot$&&$\cdot$&\\
&$\cdot$&&$\cdot$&\\
$\be(5,2,3,4)$&$=$&$\be(5,2,2,4)-\bF_{3,d_{3,u_3}}$&$=$&
	$\la 0,0,0,0;0,0;0;0,0,0\ra$
\end{tabular}
\end{center}
\caption{$\be(\bu)$ for $\be=\la 5,4,3,1;2,2;3;4,1,1\ra$
and various $\bu$}\label{examp}
\end{figure}

Finally, let $\be({\bm 0}^m)=\be$, and
denote the $k^{\rm th}$ characteristic vector $\bx_k$, having all
zeros with a single one in the $k^{\rm th}$ entry.
For ${\bm 0}^m\le\bu\le\bn$ define
$$\be(\bu)=\be(\bu-\bx_i)-\bF_{i,d_{i,u_i}}.\ $$
(Note that this definition is valid for every $1\le i\le t$.)
Figure \ref{examp} shows an example for these definitions.
Note that we always have $\be(\bn)={\bm 0}^n$.

%
%
\subsection{Lattice Graph and Pebbling Number}\label{lattice}

Define the lattice $L=L(G)=\prod_{i=1}^t P_{n_i+1}$ (the cartesian
product of paths with $n_i+1$ vertices).
Note that $L$ is isomorphic to the divisor lattice of 
$N=N(G)=\prod_{i=1}^t p_i^{e_{i,1}}$ (having height $n=\sum_{i=1}^t e_{i,1}$) 
and label the vertices of $L$ accordingly.
Next consider an edge of $L$ between vertex $p_i^kq$ and vertex
$p_i^{k-1}q$, where $p_i\not|\ q$.
Label such an edge by weight $p_i^{d_{i,k}}$.

\begin{figure}
\centerline{\psfig{figure=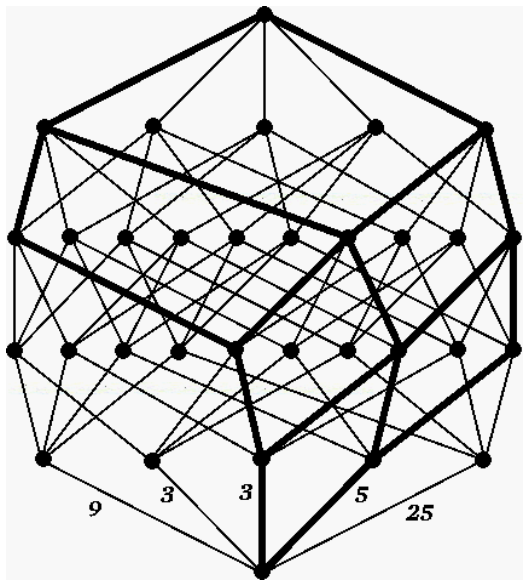,height=3.0in}}
\caption{$L(G)$ as a retract of $\B^5(9,3,3,25,5)$ for the group
$G=\aZ_9\oplus\aZ_3\oplus\aZ_3\oplus\aZ_{25}\oplus\aZ_5$}\label{retract}
\end{figure}

Because $L$ and its labelling is a retract (see Figure \ref{retract} for
an example) of the $n$-dimensional
boolean lattice $\B^n(\bw)$, having edge labels 
$\bw=\la p_i^{d_{i,j}}\ra_{i,j}$, 
we have that the generalized (pebbling operations obey the edge labels)
pebbling number $\pi(L)=\pi(\B^n(\bw))$.
(This is the same argument used in \cite{Chu}.)
Notice that
$$\pi(\B(\bw))=\prod_{i=1}^t\prod_{j=1}^{n_i}p_i^{d_{i,j}}
=\prod_{i=1}^t\prod_{j=1}^{m_i}p_i^{e_{i,j}}=|G|\ .$$ 

Given a sequence of elements of $G$, $(g_1,\ldots,g_{|G|})$,
define a configuration by placing corresponding pebbles
$\{g_1\},\ldots,\{g_{|G|}\}$ on $L$, with pebble $\{g_k\}$ on
vertex $|g_k|\in V(L)$.
Because $\pi(L)=|G|$, the configuration is solvable to the bottom 
vertex labelled 1.
As was noted in \cite{CHH}, $L$ is greedy, meaning that we may assume 
that every pebbling step moves toward the root 1.

We will now use the solution of the configuration to construct a
subsequence $(g_k)_{k\in K}$ that satisfies $\sum_{k\in K}g_k=0_G$ and
$\sum_{k\in K}1/|g_k|\le 1$.
(We will follow somewhat the structure of the argument presented in
\cite{CHH}, with a few necessary tricks thrown in.)

%
%
\subsection{Well Placed Pebbles}\label{well}

We now make several useful recursive definitions.
For a pebble $A$ define
\begin{itemize}
\item
$\set(A)=\bigcup_{B\in A}{\set}(B)$, where ${\set}(\{g_k\})=\{g_k\}$,
\item
$\val(A)=\sum_{B\in A}{\val}(B)$, where ${\val}(\{g_k\})=\bg_k$, and
\item
$\ord(A)=\sum_{B\in A}{\ord}(B)$, where ${\ord}(\{g_k\})=1/|g_k|$.
\end{itemize}
Note that ${\val}(A)=\sum_{g\in{\set}(A)}\bg$ and
${\ord}(A)=\sum_{g\in{\set}(A)}1/|g|$.
We say a pebble $A$ is {\it well placed} at vertex $\bp^{\c\ubu}$ if
\begin{enumerate}
\item
${\val}(A)\equiv {\bm 0}^m\mod \bp^{\ube(\ubu)}$ and
\item
${\ord}(A)\le 1/\bp^{\c\ubu}$.
\end{enumerate}
Thus each pebble in the initial configuration is well placed.

We will interpret each pebbling step from $x$ to $y$ as follows:
first remove a collection of pebbles $A_1,A_2,\ldots,A_s$ of the
appropriate size (the edge weight of $xy$) from $x$, then for some 
carefully chosen index set $K_x\subseteq\{1,\ldots,s\}$ place the new 
pebble $A=\{A_k\}_{k\in K_x}$ on $y$.
We will show that if each $A_k$ is well placed at $x$ then
$A$ is well placed at $y$.
Any pebble $A$ that is well placed at vertex $1=\bp^{\c{\fn\bm 0}}$ 
yields the solution ${\set}(A)$ to Theorem \ref{groups}.

%
%
\section{Proof of Theorem \ref{groups}}\label{proof}

For the purposes of notational readability, we will first give the 
proof of Theorem \ref{groups} in the case of $p$-groups.
Once established, the general case will be straightforward.

%
%
\subsection{Characteristic $p$}\label{primechar}

Here we have $t=1$ so that $i=1$ always. For ease of notation we
will simply drop the 1; thus $G=\prod_{j=1}^{m}\aZ_{p^{e_{j}}}$
for some prime $p$, multiplicity $m$, and exponents $e_{j}$
($e_1\ge\cdots\ge e_m$).  
For $\be=\la e_j\ra_{j=1}^m$ recall that
$\bp^{\c\ube}=\prod_{j=1}^mp^{e_j}=|G|$.

\begin{lem}\label{power}
Theorem \ref{groups} holds for groups of the form $G=\aZ_p^m=
\oplus_{j=1}^m\aZ_p$.
\end{lem}

\noindent
{\it Proof.}
This result will follow from Theorem \ref{cyclic}.  
View $G$ as the $m$-dimensional vector space over $\aF_p$.  
Then assign to $\aF_p^m$ the natural correspondence with field $\aF_{p^m}$,
and partition $\aF_{p^m}-\{ 0\}$ into $(p^m-1)/(p-1)$ lines of size $p-1$.

With $p^m$ pebbles, none of which is at ${\bm 0}$ (otherwise we are done), 
the pigeonhole principle forces some line to have at least $p$ pebbles.  
Since a line plus the origin forms the cycle $\aZ_p$, Theorem
\ref{cyclic} completes the proof.
\pf

\begin{thm}\label{prime}
Theorem \ref{groups} holds for groups of the form
$G=\oplus_{j=1}^{m}\aZ_{p^{e_{j}}}$.
\end{thm}

\begin{figure}
\begin{center}
$\be=\la 5,2,2,1\ra\ ,\quad\bd=\la 4,3,1,1,1\ra$
\begin{tabular}{lclcl}
\hline
&&&&\\
$\be(0)$&&&$=$&$\la 5,2,2,1\ra$\\
&&&&\\
$\be(1)$&$=$&$\la 5,2,2,1\ra-\bmf_4$&$=$&$\la 4,1,1,0\ra$\\
&&&&\\
$\be(2)$&$=$&$\la 4,1,1,0\ra-\bmf_3$&$=$&$\la 3,0,0,0\ra$\\
&&&&\\
$\be(3)$&$=$&$\la 3,0,0,0\ra-\bmf_1$&$=$&$\la 2,0,0,0\ra$\\
&&&&\\
$\be(4)$&$=$&$\la 2,0,0,0\ra-\bmf_1$&$=$&$\la 1,0,0,0\ra$\\
&&&&\\
$\be(5)$&$=$&$\la 1,0,0,0\ra-\bmf_1$&$=$&$\la 0,0,0,0\ra$\\
\end{tabular}
\end{center}
\caption{$\be(u)$ for $\be=\la 5,2,2,1\ra$
and $u=0,\ldots,5$}\label{exam2}
\end{figure}

\noindent
{\it Proof.}
We use Lemma \ref{power} to show that each pebbling step preserves
the well placed property.
Given a sequence of $|G|=\bp^{\c\ube}=\prod_{j=1}^mp^{e_j}$
elements of $G$ place, as discussed in Section \ref{well}, 
corresponding pebbles on the lattice $L=L(G)=P_{e_1+1}$, having edge 
label $p^{d_k}$ between vertices $p^k$ and $p^{k-1}$, where 
$\bd=\la d_k\ra_{k=1}^{e_1}$ is the dual partition to $\be$.
For $r\ge 0$ recall that $\bmf_r=\la {\bm 1}^r, {\bm 0}^{n-r}\ra$.
Let $\be(0)=\be$, and for $0<u\le e_1$ define
$\be(u)=\be_{u-1}-\bmf_{d_u}$ (see Figure \ref{exam2} for an example).
recall that we always have $\be(e_1)={\bm 0}^m$ because of the Ferrer's
duality.

Given $p^{d_u}$ well placed pebbles $\{A_r\}_{r=1}^{p^{d_u}}$
on vertex $p^u$, we know that each $\val(A_r)\equiv {\bm 0}^m\mod
\bp^{\ube(u)}$ and each $\ord(A_r)\le 1/p^u$.  
Consider, for each $r$, $\bB_r=\val(A_r)/\bp^{\ube(u)}$.
By Lemma \ref{power} we can find a nonempty index set $R$ so that
for $B=\{\bB_r\}_{r\in R}$ we have $\val(B)\equiv {\bm
0}^m\mod\bp^{\ubmf_{d_u}}$ and $\ord(B)\le 1$.

Now let $A=\{A_r\}_{r\in R}$.  
Then
\begin{eqnarray*}
\val(A)&=&\sum_{r\in R}\val(A_r)\\
&&\\
&=&\sum_{r\in R}\bp^{\ube(u)}\bB_r\\
&&\\
&=&\bp^{\ube(u)}\val(B)\\
&&\\
&\equiv &{\bm 0}^m\mod\bp^{\ube(u)+\ubmf_{d_u}}\\
&&\\
&=&{\bm 0}^m\mod\bp^{\ube(u-1)}\ .
\end{eqnarray*}
Also, $\ord(A)=\sum_{r\in R}\ord(A_r)\le |R|/p^u=1/p^{u-1}$.
Hence $A$ is well placed on vertex $p^{u-1}$.

Since the pebbling number guarantees that some pebble $A$ reaches
vertex $1=p^0$, and since the previous argument ensures that $A$
is well placed, we find, for some $K\not=\emt$ that
$$\sum_{k\in K}\bg_k\ =\ \val(A)\ 
\equiv\  {\bm 0}^m\mod\bp^{\ube({\fn\bm 0})}
\ =\ {\bm 0}^m\mod\bp^\be \ =\  {\bm 0}_G$$
(i.e. $\sum_{k\in K}g_k=0_G$) and 
$$\sum_{k\in K}1/|g_k|\ =\ \ord(A)\ \le \ 1/{\bp^{\c{\fn\bm 0}}}\ =\ 1\ .$$

%
%
\subsection{General Case}\label{general}

As expected, the same proof carries through; only the notation generalizes.  
Given $p^{d_{i,u_i}}_i$ well placed pebbles
$\{A_r\}_{r=1}^{p^{d_{i,u_i}}_i}$ on vertex $\bp^{\c\ubu}$, we know that
each $\val(A_r)\equiv {\bm 0}^m\mod \bp^{\ube(\ubu)}$ and each
$\ord(A_r)\le 1/\bp^{\c\ubu}$.  
Consider, for each $r$, $\bB_r=\val(A_r)/\bp^{\ube(\ubu)}$.
By Lemma \ref{power} we can find a nonempty index set $R$ so that
for $B=\{\bB_r\}_{r\in R}$ we have $\val(B)\equiv {\bm
0}^m\mod\bp^{\ubF_{i,d_{i,u_i}}}$ and $\ord(B)\le 1$.

Now let $A=\{A_r\}_{r\in R}$.  
Then
\begin{eqnarray*}
\val(A)&=&\sum_{r\in R}\val(A_r)\\
&&\\
&=&\sum_{r\in R}\bp^{\ube(\ubu)}\bB_r\\
&&\\
&=&\bp^{\ube(\ubu)}\val(B)\\
&&\\
&\equiv &{\bm 0}^m\mod\bp^{\ube(\ubu)+\ubF_{i,d_{i,u_i}}}\\
&&\\
&=&{\bm 0}^m\mod\bp^{\ube(\ubu-\ux_i)}\ .
\end{eqnarray*}
Also, $\ord(A)=\sum_{r\in R}\ord(A_r)\le
|R|/{\bp^{\c\ubu}}=1/{\bp^{\c(\ubu-\ux_i)}}$.
Hence $A$ is well placed on vertex $\bp^{\ubu-\ux_i}$.

Since the pebbling number guarantees that some pebble $A$ reaches
vertex $1=\bp^{\fn\bm 0}$, and since the previous argument ensures that $A$
is well placed, we find, for some $K\not=\emt$ that
$$\sum_{k\in K}\bg_k\ =\ \val(A)\ 
\equiv \ {\bm 0}^m\mod\bp^{\ube({\fn\bm 0})}
\ =\ {\bm 0}^m\mod\bp^\ube \ = \ {\bm 0}_G$$ 
(i.e. $\sum_{k\in K}g_k=0_G$) and
$$\sum_{k\in K}1/|g_k|\ =\ \ord(A)\ \le \ 1/{\bp^{\c{\fn\bm 0}}}\ =\ 1\ .$$

%
%
\section{Further Comments}\label{comments}

For cyclic groups Theorem \ref{groups} is best possible.
However, for other groups it is conceivable that shorter sequences of
elements may suffice.
It had been conjectured that $D(G)=1+\sum_{i=1}^r(N_i-1)$ for abelian
$G$ \cite{O}.
While this was shown true for groups of rank at most 2 and for $p$-groups,
among other special cases, it has been shown false in general
\cite{EBK,GS}.
One may ask for the generalized Davenport constant for the minimum
length of a sequence required to force a zero-sum subsequence with the
extra condition on its orders.

%
%
\section*{Acknowledgement}\label{ack}

The first author is grateful for the support of the Jack H. Hawes
Scholarship that provided him the opportunity to work on this research.

%
%
\bibliographystyle{plain}
%

%
%
\end{document}